\nonstopmode
\documentclass[12pt]{article}

\usepackage{latexsym,amssymb,amsmath}
\usepackage{amsthm, amstext}
\usepackage{array, amsfonts, mathrsfs}
\usepackage{hyperref}
\usepackage{indentfirst}
\usepackage{color}
\usepackage{ifpdf}
\usepackage{url}

\usepackage[T2A]{fontenc}
\usepackage[utf8]{inputenc}
\usepackage[russian,english]{babel}

\usepackage{amsfonts}
\usepackage{amssymb}
\usepackage{amsmath}
\usepackage{pdfpages}
\usepackage{graphicx}
\usepackage{array}
\usepackage{stackrel}
\usepackage{hhline}
\usepackage{doi}
\usepackage{caption}
\usepackage{mleftright}
\usepackage{csquotes}

\usepackage[
    backend=biber,
    clearlang=true,
    natbib=true,
    url=true,
    doi=true,
    eprint=true,
    sorting=nyt,
    autolang=hyphen,
    related=true
]{biblatex}

\newcommand{\half}{\tfrac{1}{2}}
\newcommand{\poch}[2]{\left\{\begin{array}{c}
  #1\\#2\end{array}\right\}}
  
\newcommand{\myi}{{\mathrm i}}
\newcommand{\myd}{{\mathrm d}}
\newcommand{\mye}{{\mathrm e}}
\newcommand{\abs}[1]{\vert #1 \vert}

\begin{document}
\thispagestyle{empty}
\

\vspace{-14mm}
\hfill{\bf  PDMI Preprint 07/2025}
\vskip1.5cm

\begin{center}{\Large {\bf
An approach to
the Lindel\"of hypothesis\\[3.5mm]
 for
Dirichlet $L$-functions}}
\vskip1cm

{\sc \large  Yuri Matiyasevich}
 \vskip0.7cm

{St.\,Petersburg Department\\ of
V.\,A.\,Steklov   Institute of   Mathematics
\\of    Russian   Academy   of   Sciences}

\bigskip

            {\tt{yumat@pdmi.ras.ru}}

 \vskip 2cm
 \parbox{.9\textwidth}{  \qquad\indent{\bf \small Abstract.}
  {\small
  The suggested approach is based on a known 
  representation of Dirichlet $L$-functions 
  via the incomplete  gamma functions. 
   Some properties of the Taylor coefficients of the lower incomplete gamma function at infinity seem to be new. Specifically, these coefficients can be expressed in terms of Touchard polynomials.
  Furthermore, these same coefficients can be used to reformulate the functional equation for Dirichlet $L$-functions.
  This relationship ``explains'' 
  why $\abs{L_\chi(1/2+\myi t)}$ should be small.
  
  \qquad   To present the  new ideas in a nutshell,
 we start by giving (in Section~\ref{sect1}) 
   a ``formula proof''
   of the Lindel\"of hypothesis.
This is not a genuine proof, as we are not concerned with
the convergence of our series nor do we justify changing the order of summation. 
  
  \qquad   In Section \ref{secplan},
    we suggest some  hypothetical ways of transforming the ``proof'' from Section \ref{sect1} into a rigorous mathematical proof.
 
   \qquad   Sections 3--5  contain some 
   technical detail and bibliographical references.

}}

\end{center}

\vskip1cm


 \newpage
 
\thispagestyle{empty}
{\centering
{\bf ПРЕПРИНТЫ ПОМИ РАН}\\[7mm]

{\bf PREPRINTS
OF THE ST. PETERSBURG DEPARTMENT\\ OF STEKLOV INSTITUTE OF MATHEMATICS}\\[12mm]

ГЛАВНЫЙ РЕДАКТОР\\[3mm]
{\sc С.\,В.\,Кисляков}\\[13mm]
РЕДКОЛЛЕГИЯ\\[3mm]\small
{\sc В.\,М.\,Бабич,   М.\,А.\,Всемирнов, А.\,И.\,Генералов,
И.\,А.\,Ибрагимов, Л.\,Ю.\,Колотилина,  Ю.\,В.\,Матиясевич,
Н.\,Ю.\,Нецветаев, С.\,И.\,Репин, Г.\,А.\,Серегин}\\[13mm]
Учредитель: Федеральное государственное бюджетное
учреждение науки

Санкт-Петербургское отделение
Математического института

им.\ В.\,А.\,Стеклова Российской академии наук\\[13mm]
Свидетельство о регистрации средства массовой информации:

 ЭЛ №ФС 77-33560 от 16 октября 2008 г.

  Выдано Федеральной службой по надзору\\в сфере связи и массовых коммуникаций\\[13mm]
Контактные данные:\\[2mm] 191023, г.\ Санкт-Петербург, наб.\ реки Фонтанки, дом 27\\[2mm]
телефоны:\ (812) 312-40-58; (812) 571-57-54\\[2mm]
e-mail: {\tt admin@pdmi.ras.ru}\\[2mm]
{\tt https://www.pdmi.ras.ru/preprint/}\\[8mm]
Заведующая информационно-издательским сектором {\sc В.\,Н.\,Симонова}

}

\newpage

\vspace{30mm}

\section{``Formula proof''}\label{sect1}

We start by presenting  what could be called 
   a 'formula proof'
   of the Lindel\"of hypothesis.
In this section, we ignore 
the convergence of our series and  do we  justify changing the order of summation. In the next section, we discuss possible ways to
transform this ``proof'' into a rigorous mathematical argument.

To simplify notation, this section deals with a 
specific
function $L_3(s)$; a generalisation to  
an arbitrary Dirichlet $L$-function
is considered in Section \ref{thetagen}.

The function $L_3(s)$ is defined (for $\Re(s)>0$)
by the Dirichlet series.
\begin{equation}\label{L3}
  L_3(s)=\sum_{n=1}^\infty \chi_3(n)n^{-s},
\end{equation}
where $\chi_3$ is the only non-principal character 
modulo~$3$:
\begin{equation}\label{chi3}
 \chi_3(n)=
 \begin{cases}
   \phantom{-}0,& \text{if }n\equiv 0 \pmod{3},\\
   \phantom{-}1,& \text{if } n\equiv 1 \pmod{3},\\
   -1,& \text{if } n\equiv 2 \pmod{3}.
\end{cases}
\end{equation}

\

\textbf{The Lindel\"of hypothesis (for the case of $L_3(s)$).}
\emph{For every positive  $\varepsilon$
\begin{equation}\label{Lindelof}
  L_3(\tfrac{1}{2}+\myi t)=O_\varepsilon(t^\varepsilon)
\end{equation}
as $t\rightarrow+\infty$.}

\  

The hypothesis can also be stated in terms of function.
\begin{equation}\label{xigL}
  \xi_3(s)=g_3(s)L_3(s),
\end{equation}
where
\begin{equation}\label{g}
  g_3(s)=\left(\tfrac{\pi}{3}\right)^{-\frac{s+1}{2}}
  \Gamma\mleft(\tfrac{s+1}{2}\mright).
\end{equation}
Namely, for a positive $a$  
\begin{equation}\label{growthgamma}
  \vert\Gamma(a+\myi\,t)\vert=(1+o_a(t))
  \sqrt{2 \pi }
  t^{a-\frac{1}{2}}
  \mye^{-\frac{ \pi   t} {2}}\quad \text{as }
   t\rightarrow+\infty.
\end{equation}
Thus \eqref{Lindelof} is equivalent to
 \begin{equation}\label{Lindelofxi}
  \xi_3(\tfrac{1}{2}+\myi t)=O_\epsilon\mleft(
  t^{\frac{1}{4}+\epsilon}
  \mye^{-\frac{\pi    t}{4}}\mright).
\end{equation}

It is well known that the function $\xi_3(s)$ satisfies the
\emph{functional equation}
\begin{equation}\label{fe}
  \xi_3(s)=\xi_3(1-s).
\end{equation}
This implies that 
\begin{equation}\label{xixixi}
  \xi_3(\tfrac{1}{2}+\myi t)=
  \half \xi_3(\tfrac{1}{2}+\myi t)+\half \xi_3(\tfrac{1}{2}-\myi t).
\end{equation}
Straightforward substitution of 
\eqref{L3} and \eqref{g} into \eqref{xixixi}
gives the following:
\begin{multline}\label{splitted}
\xi_3(\tfrac{1}{2}+\myi t)=
\frac{1}{2}\sum_{n=1}^\infty\chi_3(n)
    { \left(\tfrac{\pi}{3}\right)}^{-\frac{3}{4}-\frac{\myi t}{2}}
  \Gamma
   \left(\tfrac{3}{4}+\tfrac{\myi t}{2}\right)
   n^{-\frac{1}{2}-\myi t}+\\
   +
   \frac{1}{2}\sum_{n=1}^\infty\chi_3(n)
    { \left(\tfrac{\pi}{3}\right)}^{-\frac{3}{4}+\frac{\myi t}{2}}
  \Gamma
   \left(\tfrac{3}{4}-\tfrac{\myi t}{2} 
   \right)
   n^{-\frac{1}{2}+\myi t}.
\end{multline}
It turns out (for details see Section \ref{Lavrik})
 that we can drop the  factors $\half$ and replace the
complete gamma function with the lower incomplete gamma function
with a suitable second argument:
\begin{multline}\label{splitted2}
\xi_3(\tfrac{1}{2}+\myi t)=
\sum_{n=1}^\infty\chi_3(n)
    { \left(\tfrac{\pi}{3}\right)}^{-\frac{3}{4}-\frac{\myi t}{2}}
  \gamma
   \left(\tfrac{3}{4}+\tfrac{\myi t}{2},
   \tfrac{ \pi   n^2  }{3}\right)
   n^{-\frac{1}{2}-\myi t}+\\
   +
   \sum_{n=1}^\infty\chi_3(n)
    { \left(\tfrac{\pi}{3}\right)}^{-\frac{3}{4}+\frac{\myi t}{2}}
  \gamma
   \left(\tfrac{3}{4}-\tfrac{\myi t}{2} , \tfrac{ \pi   n^2  }{3}\right)
   n^{-\frac{1}{2}+\myi t}.
\end{multline}

The following expansion (due to E. E. Kummer) holds for the lower gamma function:
\begin{equation}\label{kumm}
  \gamma(w,m)=m^w\sum_{k=0}^\infty
  \frac{(-1)^k}{k!(w+k)}m^k.
\end{equation}
Assuming that
\begin{equation}\label{geo}
  \frac{1}{w+k}=\sum_{l=1}^\infty
  (-1)^{l-1}\frac{k^{l-1}}{w^l},
\end{equation}
we have
\begin{eqnarray}\label{sum1}
   \gamma(w,m)
   &=&
   m^w\sum_{k=0}^\infty
    \sum_{l=1}^\infty
  \frac{(-1)^{k+l-1}k^{l-1}}{k!w^l}
    m^k\\\label{sum2}
 &=&
   m^w
    \sum_{l=1}^\infty
    \frac{(-1)^{l-1}}{w^l}
    \sum_{k=0}^\infty
  \frac{(-1)^{k}k^{l-1}}{k!}
    m^k
\end{eqnarray}
(in \eqref{geo}--\eqref{sum2} and in the sequel
we assume that $0^0=1$).
Further we have:
\begin{equation}\label{power}
  k^l=\sum_{j=0}^l\poch{l}{j}k^{\underline{j}},
\end{equation}
where $\poch{l}{j}$ is 
the \emph{Stirling number of the second kind} 
and  $k^{\underline{j}}$ is the \emph{falling 
factorial power},
\begin{equation}\label{fall}
  k^{\underline{j}}=k(k-1)\dots(k-j+1)
\end{equation}
(we use the notation from \cite{KnuthGP}).

Substituting
 \eqref{power} into \eqref{sum2}, we get:
\begin{eqnarray}
   \gamma(w,m)
   &=&
      \label{sum4}
     m^w
    \sum_{l=1}^\infty
    \frac{(-1)^{l-1}}{w^l}
    \sum_{k=0}^\infty
  \frac{(-1)^{k}}{k!}
    m^k\sum_{j=0}^{l-1}\poch{l-1}{j}k^{\underline{j}}
          \\  \label{sum5} &=&
     m^w
    \sum_{l=1}^\infty
    \frac{(-1)^{l-1}}{w^l}
    \sum_{j=0}^{l-1}\poch{l-1}{j}
    \sum_{k=0}^\infty
    (-1)^{k}
  \frac{k^{\underline{j}}}{k!}
    m^k
             \\ \label{sum6}  &=&
     m^w
    \sum_{l=1}^\infty
    \frac{(-1)^{l-1}}{w^l}
    \sum_{j=0}^{l-1}(-1)^j\poch{l-1}{j}
    m^j\sum_{k=j}^\infty
  \frac{ (-1)^{k-j}}{(k-j)!}
    m^{k-j}
    \\ \label{sum8}  &=&
     m^w \mye^{-m}
    \sum_{l=1}^\infty
    \frac{(-1)^{l-1}T_{l-1}(-m)}{w^l},
   \end{eqnarray}
where
\begin{equation}\label{touchard}
  T_n(x)=  \sum_{j=0}^n\poch{n}{j}x^j.
\end{equation}

Polynomials \eqref{touchard} are known as 
 \emph{Touchard polynomials} (after \cite{Touchard1939}). 
They have been widely studied, but the  author 
has not found the equality~\eqref{sum8} in the literature. 
This equality can be viewed from different points:
\begin{itemize}
  \item as
the Taylor expansion of the lower gamma function 
at infinity;
\item 
  as a definition of Touchard polynomials 
with $m^{-w}\mye^m\gamma(w,m)$ playing  the role of
a generating function. 
\end{itemize}

\

Substituting  $a+\myi b$ for $w$ in \eqref{sum8}, we get that
\begin{eqnarray}
  \label{sum9}
   \gamma(a+\myi b,m)&=&
     m^{a+\myi b} \mye^{-m}
    \sum_{l=1}^\infty
    \frac{(-1)^{l-1}T_{l-1}(-m)}{(a+\myi b)^l}.
   \end{eqnarray}
Assuming that
\begin{equation}\label{geol}
  \frac{1}{(a+\myi b)^l}=
  {(-\myi)^l}
  \sum_{j=0}^\infty
  \binom{l+j-1}{j}\frac{(\myi a)^{j}}{b^{l+j}},
\end{equation}
we get that
\begin{eqnarray}
  \label{sum10}
   \gamma(a+\myi b,m)&=&
     -m^{a+\myi b} \mye^{-m}
    \sum_{l=1}^\infty
    {\myi^{l}
   T_{l-1}(-m)}
    \sum_{j=0}^\infty
   \binom{l+j-1}{j}\frac{(\myi a)^{j}}{b^{l+j}}\\
   \label{sum11}&=&
   m^{a+\myi b} \mye^{-m}
     \sum_{k=1}^\infty
    \frac{U_k(a,m)}{b^k},
   \end{eqnarray}
where
\begin{equation}\label{defU}
  U_k(a,m)=-{\myi^{k}
  \sum_{j=0}^{k-1}
   \binom{k-1}{j}{ a^{j}}}
    T_{k-j-1}(-m).
\end{equation}

   Let us set $a=\tfrac{3}{4}$ and  $b=\pm\tfrac{ t}{2}$
    in \eqref{sum11}. Substituting the 
    resulting right-hand sides of  \eqref{sum11} into 
    \eqref{splitted2}, we get that
\begin{eqnarray}\label{xiU1}
  \xi_3(\tfrac{1}{2}+\myi t)
  &=&
  \sum_{n=1}^\infty \chi_3(n)n \mye^{-\tfrac{ \pi   n^2  }{3}}
 \sum_{k=1}^\infty {2^k U_{k}(\tfrac{3}{4},
  \tfrac{ \pi   n^2  }{3})}
  \left(\tfrac{1}{t^k}+\tfrac{1}{(-t)^k}\right)\\
  \label{xiU12}
  &=&
  \sum_{n=1}^\infty \chi_3(n)n \mye^{-\tfrac{ \pi   n^2  }{3}}
 \sum_{k=1}^\infty { U_{2k}(\tfrac{3}{4},
  \tfrac{ \pi   n^2  }{3})}
  \tfrac{2^{2k+1}}{t^{2k}}\\
  \label{xiU2}
  &=&\sum_{k=1}^\infty
  \frac{2^{2k+1}}{t^{2k}}
  \sum_{n=1}^\infty \chi_3(n)n
 {U_{2k}(\tfrac{3}{4},
  \tfrac{ \pi   n^2  }{3})\mye^{-\tfrac{ \pi   n^2  }{3}}}.
\end{eqnarray}

Our approach to the Lindel\"of hypothesis 
is based on the following 

\

\textbf{Key discovery.}
\emph{For every $k$}
\begin{eqnarray}\label{key}
\sum_{n=1}^\infty \chi_3(n)n
 {U_{2k}(\tfrac{3}{4},
  \tfrac{ \pi   n^2  }{3})\mye^{-\tfrac{ \pi   n^2  }{3}}}
  =0.
\end{eqnarray}

\

The equalities \eqref{key} are 
essentially equivalent to the functional
equation \eqref{fe} (see Section \ref{thetapol}).

\

\textbf{Corollary of \eqref{xiU2} and \eqref{key}.} \emph{For every real $t$}
\begin{equation}
\label{xiis0}
  \xi_3(\tfrac{1}{2}+\myi t)=0.
\end{equation}

\

The  desired equality 
\eqref{Lindelofxi} follows from the  paradoxical identity  \eqref{xiis0}.
However, the latter is, of course, incorrect. Nevertheless,   the above
formal  deduction of   \eqref{xiis0}  suggests 
an approach to the Lindel\"{o}f hypothesis (see the next section).

\section{Approaches to the  Lindel\"{o}f hypothesis}
\label{secplan}

Here we consider hypothetical ways in which considerations from
the previous section could be transformed into rigorous
estimates of $\abs{ \xi_3(\tfrac{1}{2}+\myi t)}$.

Our first error was \eqref{geo}. The equality holds when 
$k<\abs{w}$ only. Later, we did not justify the change of order of summation in
\eqref{sum1}--\eqref{sum2},
\eqref{sum4}--\eqref{sum5},
 \eqref{sum10}--\eqref{sum11}, and
 \eqref{xiU1}--\eqref{xiU2}.
Nevertheless, \eqref{sum8} and \eqref{sum11} are
 asymptotic representations 
for the lower gamma function, that is,
\begin{eqnarray}\label{sum8a}
   \gamma(w,m)
   &=&
     m^w \mye^{-m}\left(
    \sum_{l=1}^L
    \frac{(-1)^{l-1}T_{l-1}(-m)}{w^l}
    +O_{m,L}\left(\tfrac{1}{w^{L+1}}\right)\right)
   \end{eqnarray}
   and
 \begin{eqnarray}
     \gamma(a+\myi b,m)
    &=&\label{sum11a}
   m^{a+\myi b} \mye^{-m}\left(
     \sum_{k=1}^K
    \frac{U_k(a,m)}{b^k}+O_{a,m,K}\left(\tfrac{1}{b^{K+1}}\right)\right).
   \end{eqnarray}

As a natural attempt to improve our ``formula proof'' we could
introduce two bounds, $T$ and $K$
(depending on $t$, $N=N(t)$, $K=K(t)$),  substitute them 
for the  upper 
limits (equal to the infinity) in the summations in \eqref{splitted2}
and  \eqref{xiU1}--\eqref{xiU2}, and try to estimate  the emerging errors. 

Let (cf.\ \eqref{splitted2}) 
\begin{multline}\label{splitted2N}
\xi_{3,N}(\tfrac{1}{2}+\myi t)=
\sum_{n=1}^N\chi_3(n)
    { \left(\tfrac{\pi}{3}\right)}^{-\frac{3}{4}-\frac{\myi t}{2}}
  \gamma
   \left(\tfrac{3}{4}+\tfrac{\myi t}{2},
   \tfrac{ \pi   n^2  }{3}\right)
   n^{-\frac{1}{2}-\myi t}+\\
   +
   \sum_{n=1}^N\chi_3(n)
    { \left(\tfrac{\pi}{3}\right)}^{-\frac{3}{4}+\frac{\myi t}{2}}
  \gamma
   \left(\tfrac{3}{4}-\tfrac{\myi t}{2} , \tfrac{ \pi   n^2  }{3}\right)
   n^{-\frac{1}{2}+\myi t}
\end{multline}
and (cf.\ \eqref{xiU12}--\eqref{xiU2}) 
\begin{eqnarray}\label{xiU1NK}
  \xi_{3,N,K}(\tfrac{1}{2}+\myi t)&=&
  \sum_{n=1}^N \chi_3(n)n \mye^{-\tfrac{ \pi   n^2  }{3}}
 \sum_{k=1}^K { U_{2k}(\tfrac{3}{4},
  \tfrac{ \pi   n^2  }{3})}
  \tfrac{2^{2k+1}}{t^{2k}}\\
  \label{xiU2KN}
  &=&\sum_{k=1}^K
  \tfrac{2^{2k+1}}{t^{2k}}
  \sum_{n=1}^N \chi_3(n)n
 {U_{2k}(\tfrac{3}{4},
  \tfrac{ \pi   n^2  }{3})\mye^{-\tfrac{ \pi   n^2  }{3}}}.
\end{eqnarray}
We have:
\begin{eqnarray}\label{diff0}
  \lefteqn{\hspace{7mm}\abs{\xi_{3}(\tfrac{1}{2}+\myi t)}\le}\nonumber\\\label{diff1}
  &&\abs{\xi_{3}(\tfrac{1}{2}+\myi t)-
  \xi_{3,N}(\tfrac{1}{2}+\myi t)}+\\\label{diff2}
   &&\hspace{18mm} +\abs{\xi_{3,N}(\tfrac{1}{2}+\myi t)-
  \xi_{3,N,K}(\tfrac{1}{2}+\myi t)}+\\
   &&\hspace{43mm}+
  \abs{\xi_{3,N,K}(\tfrac{1}{2}+\myi t)}.\label{diff3}
\end{eqnarray}

 Thanks to the convergence of the series in \eqref{splitted2},
 the value of
 \eqref{diff1}
 can  be made 
arbitrarily small.
 However, this requires the selection of a sufficiently 
 large $N$.

 Thanks to \eqref{key}, 
 each inner sum in \eqref{xiU2KN}
 tends to zero as $N\rightarrow\infty$.
Hence  for a fixed $K$ the value of
\eqref{diff3}
can be made 
arbitrarily small, again at the cost of 
selecting sufficiently large $N$.  

The choice of $K$ is trickier. We have:
\begin{multline}\label{xiU1NKacc}
  \xi_{3,N,K}(\tfrac{1}{2}+\myi t)=
  \sum_{n=1}^N \chi_3(n)n \mye^{-\tfrac{ \pi   n^2  }{3}}
 \sum_{k=1}^{2K} { U_{k}(\tfrac{3}{4},
  \tfrac{ \pi   n^2  }{3})}
  \tfrac{2^{2k}}{t^{k}}+\\+
   \sum_{n=1}^N \chi_3(n)n \mye^{-\tfrac{ \pi   n^2  }{3}}
 \sum_{k=1}^{2K} { U_{k}(\tfrac{3}{4},
  \tfrac{ \pi   n^2  }{3})}
  \tfrac{2^{2k}}{(-t)^{k}}  ,
        \end{multline}
thus 
\begin{eqnarray}
 \lefteqn{  \xi_{3,N}
 \left(
 \tfrac{1}{2}+\myi t\right)-
   \xi_{3,N,K}
 \left(
 \tfrac{1}{2}+\myi t\right)
 =\sum_{n=1}^N \chi_3(n)\times}\nonumber\\&\times \left(
   {\pi^{-\left(\frac{3}{4}+\frac{ \myi  t}{2}\right)}}  
   \gamma\mleft(
      \tfrac{3}{4}+\tfrac{\myi t}{2},
       \tfrac{\pi n^2}{3}
   \mright)
  n^{-(\frac{1}{2}+\myi t)}   
 -
 \mye^{-\frac{\pi n^2}{3}}   
  \sum_{k=1}^{2K}
   \frac{2^k n U_k\mleft(\tfrac{3}{4}, \frac{\pi n^2}{3}\mright)}
  { t^k}\right.
  +\nonumber\\&\  \ +
  \left.
   {\pi^{-\left(\frac{3}{4}-\frac{ \myi  t}{2}\right)}}  
   \gamma\mleft(
      \tfrac{3}{4}-\tfrac{\myi t}{2},
       \tfrac{\pi n^2}{3}
   \mright)
  n^{-(\frac{1}{2}-\myi t)}   
 -
 \mye^{-\frac{\pi n^2}{3}}   
  \sum_{k=1}^{2K}
   \frac{2^k n U_k\mleft(\tfrac{3}{4}, \frac{\pi n^2}{3}\mright)}
  {(- t)^k}\right).\label{xiU1Nbis}
  \end{eqnarray}

According to \eqref{sum11a} (with  $a=\tfrac{3}{4}$,
$b=\pm\tfrac{ t}{2}$, and $m=\tfrac{\pi n^2}{3}$) 
\begin{equation} \label{Kbound}
\  \hspace{-2mm} {\pi^{-\left(\frac{3}{4}\pm\frac{ \myi  t}{2}\right)}}  
   \gamma\mleft(
      \tfrac{3}{4}\pm\tfrac{\myi t}{2},
       \tfrac{\pi n^2}{3}
   \mright)
  n^{-(\frac{1}{2}\pm\myi t)}   
 -
 \mye^{-\frac{\pi n^2}{3}}   
  \sum_{k=1}^K
   \frac{2^k n U_k\mleft(\frac{3}{4}, \frac{\pi n^2}{3}\mright)}
  {(\pm t)^k}=O_{n,K}\left(\frac{1}{t^{K+1}}\right).
\end{equation} 
 Using the explicit bounds (hidden in $O_{n,K}$ notation in
 \eqref{Kbound}), we can estimate~\eqref{diff2}. However, it 
 is unclear whether this (together with the estimations of
 \eqref{diff1} and~\eqref{diff3}) would result in 
 the desired estimate \eqref{Lindelofxi}.
 
The following approach looks more promising.

There are numerous known  methods for accelerating slowly
convergent series or even obtaining 'correct' value of
a divergent series. Many of these methods are linear:
an infinite sum
\begin{equation}
  \sum_{n=1}^\infty a_n
\end{equation} 
is very well approximated 
for relatively small $N$ by a finite sum
\begin{equation}
  \sum_{n=1}^N \nu_n(N)a_n
\end{equation} 
for suitable weights $\nu_n(N)$. We can 
select two  sequences of weights, 
 $\nu_n(N)$ and  $\kappa_k(K)$,  
and redefine \eqref{splitted2N}, \eqref{xiU1NK},
and \eqref{xiU2KN} as weighted finite sums:
\begin{multline}\label{splitted2Nmu}
\xi_{3,N}(\tfrac{1}{2}+\myi t)=
\sum_{n=1}^N\nu_n(N)\chi_3(n)
    { \left(\tfrac{\pi}{3}\right)}^{-\frac{3}{4}-\frac{\myi t}{2}}
  \gamma
   \left(\tfrac{3}{4}+\tfrac{\myi t}{2},
   \tfrac{ \pi   n^2  }{3}\right)
   n^{-\frac{1}{2}-\myi t}+\\
   +
   \sum_{n=1}^N\nu_n(N)\chi_3(n)
    { \left(\tfrac{\pi}{3}\right)}^{-\frac{3}{4}+\frac{\myi t}{2}}
  \gamma
   \left(\tfrac{3}{4}-\tfrac{\myi t}{2} , \tfrac{ \pi   n^2  }{3}\right)
   n^{-\frac{1}{2}+\myi t},
\end{multline}
\begin{eqnarray}\label{xiU1NKbis}
  \xi_{3,N,K}(\tfrac{1}{2}+\myi t)&=&
  \sum_{n=1}^N \nu_n(N)\chi_3(n)n \mye^{-\tfrac{ \pi   n^2  }{3}}
 \sum_{k=1}^K \kappa_k(K){ U_{2k}(\tfrac{3}{4},
  \tfrac{ \pi   n^2  }{3})}
  \tfrac{2^{2k+1}}{t^{2k}}\\
  \label{xiU2KNbis}
  &=&\sum_{k=1}^K\kappa_k(K)
  \tfrac{2^{2k+1}}{t^{2k}}
  \sum_{n=1}^N \nu_n(N)\chi_3(n)n
 {U_{2k}(\tfrac{3}{4},
  \tfrac{ \pi   n^2  }{3})\mye^{-\tfrac{ \pi   n^2  }{3}}}.
\end{eqnarray}

The weights $\kappa_k(K)$ should provide a good approximation
of the lower gamma function via a weighted finite sum (cf. \eqref{Kbound}):
\begin{equation} \label{Kboundacc}
\  \hspace{-2mm} {\pi^{-\left(\frac{3}{4}\pm\frac{ \myi  t}{2}\right)}}  
   \gamma\mleft(
      \tfrac{3}{4}\pm\tfrac{\myi t}{2},
       \tfrac{\pi n^2}{3}
   \mright)
  n^{-(\frac{1}{2}\pm\myi t)}   
 \approx
 \mye^{-\frac{\pi n^2}{3}}   
  \sum_{k=1}^K
  \kappa_k(K)
   \frac{2^k n U_k\mleft(\frac{3}{4}, \frac{\pi n^2}{3}\mright)}
  {(\pm t)^k}.
\end{equation}

The weights $\nu_n(N)$ should satisfy two conditions. Firstly,
the finite sums from \eqref{splitted2Nmu}
should provide a good approximation to
$\xi_3(\frac{1}{2}+\myi t)$ (cf. \eqref{splitted2}):
\begin{equation}\label{splitted2Nmuacc}
\xi_{3}(\tfrac{1}{2}+\myi t)\approx
\xi_{3,N}(\tfrac{1}{2}+\myi t)
.
\end{equation}
Secondly, the inner sums in \eqref{xiU2KNbis}
should be small (cf. \eqref{key}):
\begin{equation}
\sum_{n=1}^N \nu_n(N)\chi_3(n)n
 {U_{2k}(\tfrac{3}{4},
  \tfrac{ \pi   n^2  }{3})\mye^{-\tfrac{ \pi   n^2  }{3}}}
  \approx 0.
  \end{equation}

One can hope  that
for a suitable choice of the weights
it will be possible to use  values of 
$N$ and $K$ that are not too large, and deduce \eqref{Lindelofxi} from 
\eqref{diff1}--\eqref{diff3}
(after the redefinitions \eqref{splitted2Nmu}--\eqref{xiU2KNbis}). 

There are efficient non-linear methods for accelerating convergence. However, linearity
is essential for us: we need both representations  
\eqref{xiU1NKbis}
 and \eqref{xiU2KNbis}.

\section{Representations of $L$-functions via
incomplete gamma functions}\label{Lavrik}

A predecessor of \eqref{splitted2} can already  be found 
in the seminal paper of B.\ Riemann~\cite{Riemann1859}.
Naturally, it provides an expression 
 for the zeta function.
For Dirichlet $L$-functions, A.\,F.\,Lavrik 
obtained the following general representation.

Let $\chi(n)$ be a primitive character 
modulo  some~$q$, $q>1$,
\begin{equation}\label{xigLL}
  \xi_\chi(s)=g_\chi(s)L_\chi(s),
\end{equation}
where
\begin{equation}\label{LL}
  L_\chi(s)=\sum_{n=1}^\infty \chi(n)n^{-s},
\end{equation}
\begin{equation}\label{gL}
  g_\chi(s)=\left(\frac{\pi}{q}\right)^{-\frac{s+\delta}{2}}
  \Gamma\mleft(\frac{s+\delta}{2}\mright),
\end{equation}
\begin{equation}\label{kappa}
  \delta=
  \begin{cases}
  0,&\text{if }\chi(-1)=1,\\
  1,&\text{if }\chi(-1)=-1.
  \end{cases}
\end{equation}

The function  $\xi_\chi(s)$ satisfies the  functional equation
\begin{equation}\label{feL}
  \xi_\chi(s)=\omega\xi_{\bar{\chi}}(1-s),
  \end{equation}
where
\begin{equation}
  \bar{\chi}(n)=\overline{\chi(n)}
\end{equation}
and
\begin{equation}\label{omega}
  \omega=\frac{\sum_{k=1}^q\chi(k)
  \mye^\frac{2\pi\myi k}{q}}{\myi^\delta\sqrt{q}}.
\end{equation}  

A.\,F.\,Lavrik proved in \cite[Th.\,1]{Lavrik1966}
that for all $s$
\begin{multline}\label{splitted2L}
\xi_\chi(s)=
\sum_{n=1}^\infty\chi(n)
    { \left(\tfrac{\pi}{q}\right)}^{-\frac{s+\delta}{2}}
  \Gamma
   \left(\tfrac{s+\delta}{2},
   \tfrac{ \pi   n^2 \tau }{q}\right)
   n^{-s}+\\
   +\omega
   \sum_{n=1}^\infty\bar{\chi}(n)
    { \left(\tfrac{\pi}{q}\right)}^{-\frac{1-s+\delta}{2}}
  \Gamma
   \left(\tfrac{1-s+\delta}{2} , \tfrac{ \pi   n^2  }{q\tau}\right)
   n^{-(1-s)},
\end{multline}
where $\tau$ is an arbitrary complex number with
positive real part.

From \eqref{splitted}
and the well-known identity
\begin{equation}
  \Gamma(a)=\gamma(a,z)+\Gamma(a,z)
\end{equation}
we obtain a dual representation in which the upper gamma function is replaced by the 
lower gamma function:
 \begin{multline}\label{splitted2Low}
\xi_\chi(s)=
\sum_{n=1}^\infty\chi(n)
    { \left(\tfrac{\pi}{q}\right)}^{-\frac{s+\delta}{2}}
  \gamma
   \left(\tfrac{s+\delta}{2},
   \tfrac{ \pi   n^2 \tau }{q}\right)
   n^{-s}+\\
   +\omega
   \sum_{n=1}^\infty\bar{\chi}(n)
    { \left(\tfrac{\pi}{q}\right)}^{-\frac{1-s+\delta}{2}}
  \gamma
   \left(\tfrac{1-s+\delta}{2} , \tfrac{ \pi   n^2  }{q\tau}\right)
   n^{-(1-s)}.
\end{multline}
Representation \eqref{splitted2} is a special case of 
\eqref{splitted2Low} with 
$q=3$, $\chi=\chi_3$ (defined by \eqref{chi3}),
$\delta=1$, $\bar{\chi}=\chi_3$, $\omega=1$ and~$\tau=1$.

From a computational point of view 
representation \eqref{splitted2L}
is more efficient than \eqref{splitted2Low}.
This is because the upper gamma decays 
quickly  as the imaginary part of the second argument
tends to infinity; in this case, the lower gamma function 
tends to a non-zero limit.  However, for our purposes
the case of the lower gamma function is
more interesting because the
polynomials $U_{k}(a,m)$ from 
representation \eqref{sum11}
satisfy the identity~\eqref{key}. 

The representation \eqref{splitted2L} was generalized in a number of
ways. 
A.\,F.\,Lavrik~\cite{Lavrik1992} 
\nocite{Lavrik1992T}
and
M.\ Alzergani
\cite{Alzergani2023} gave (different) representations
which involve a finite number of factors from the Euler product.
 For another generalisation see, for example, the survey by  
 M.\ O.\ Rubinstein
 \cite{Rubinstein2004}.

\section{Corollaries of the functional equation
for~$\xi_3$}\label{thetapol}

The functional equation \eqref{fe}
can be proved in many ways. A standard technique
(see, for example, \nocite{VoroninKT}
\cite[Ch.\,1, Sect. 4, Th.\,1]{kar}) is based on the properties of Dirichlet $\theta$-functions.
Let us define $\theta_3(\tau)$
for $\Re(\tau)>0$ as follows:
\begin{equation}\label{thetaD}
  \theta_3(\tau)=\sum_{n=-\infty}^\infty
  \chi_3(n)n\mye^{-\frac{\ \pi n  ^2}{3}\tau}.
\end{equation}
This function satisfies (\cite[Ch.\,1, Sect. 4, Lemma\,2]{kar}) the functional equation
\begin{equation}\label{thetaFE}
 \theta_3({\tau}^{-1})=
\tau^{\frac{3}{2}}
  \theta_{{3}}(\tau).
\end{equation}
The functional equation \eqref{fe}
can be deduced from \eqref{thetaFE}.
A.\,F.\,Lavrik shows in \cite{Lavrik1991}\nocite{Lavrik1991T} that,
\emph{vice versa}, \eqref{thetaFE} can be
deduced from \eqref{fe}.

We  reformulate 
\eqref{thetaFE} in the following way.
This \emph{functional} equality (an identity in $\tau$)
implies 
\emph{numerical} equalities
\begin{equation}\label{numeq}
  \left.\frac{\myd^m}{\myd \tau^m} 
  \left(\theta_3({\tau}^{-1})-
{\tau^{\frac{3}{2}}}
  \theta_{{3}}(\tau)\right)\right\vert_{\tau=1}
  =0,\qquad m=1,2,\dots \ .
\end{equation}
In their turn, the numerical equalities \eqref{numeq},
taken together,  imply the functional equality~\eqref{thetaFE}.

The left-hand side in \eqref{numeq} can be written as
\begin{equation}\label{defG}
  \left.\frac{\myd^m}{\myd \tau^m} 
  \left(\theta_3({\tau}^{-1})-
{\tau^{\frac{3}{2}}}
  \theta_{{3}}(\tau)\right)\right\vert_{\tau=1}
  =\sum_{n=1}^\infty\chi_3(n)
  G_m(n)\mye^{-\frac{\ \pi n  ^2}{3}}
\end{equation}
for certain polynomials $G_m(n)$. 
Some of these are presented in Table \ref{G}. 
Explicit representations for $G_m(n)$ are given by
\eqref{GEF}, 
 \eqref{defE2} and \eqref{defF2}.

In general, the polynomial $G_m(n)$ contains only odd powers of $n$.
The degree of $G_m(n)$ is equal to $2m+1$ for an odd $m$ and
$2m-1$ for an even $m$.
 
 With the new notation, the equalities \eqref{numeq}
 can be written as
 \begin{equation}\label{sumGis0}
 \sum_{n=1}^\infty\chi_3(n)
  G_m(n)\mye^{-\frac{\ \pi n  ^2}{3}}=0.
\end{equation}
 These equalities   are not independent. For example,
$G_2(n)=2G_1(n)$.

Let ${\cal{G}}_{3}(m) $ be the linear 
   span of polynomials
 $G_1(n),\dots,G_m(n)$. For an even~$m$
the dimension of the space ${\cal{G}}_{3}(m) $  is  equal to
 $m/2$ and
 \begin{equation}\label{GG}
   {\cal{G}}_{3}(m) ={\cal{G}}_{3}(m-1) .
 \end{equation} 
 
A basis of  ${\cal{G}}_{3}(m) $ can be selected 
in many ways. For example,
for an even $m$  
polynomials
\begin{equation}\label{basisO}
 G_1(n),G_3,\dots,G_{m-1} (n)
\end{equation}
and
polynomials
\begin{equation}\label{basisE}
 G_2(n),G_4,\dots,G_{m} (n)
\end{equation}
span \eqref{GG}.

Clearly, if $P(n)\in  {\cal{G}}_{3}(m) $, then
\begin{eqnarray}\label{Pis0}
\sum_{n=1}^\infty \chi_3(n)
 {P(n)\mye^{-\tfrac{ \pi   n^2  }{3}}}
  =0.
\end{eqnarray}
Identity \eqref{key} is a corollary of

\

\textbf{Key discovery (stronger form).}
\emph{For every $k$}
\begin{eqnarray}\label{keystrong}
 nU_{2k}(\tfrac{3}{4},  \tfrac{ \pi   n^2  }{3})\in {\cal{G}}_{3}(2k).
\end{eqnarray}

\

Tables \ref{UGO} and  \ref{UGE} present a few of 
the corresponding linear relations
for bases~\eqref{basisO} and \eqref{basisE} respectively.

The relations in Tables \ref{UGO} and  \ref{UGE}  can be inverted; 
Table \ref{GU} presents a few of the inverse relations.

In the general case, the polynomials
\begin{equation}\label{basisU}
  nU_{2}(\tfrac{3}{4},  \tfrac{ \pi   n^2  }{3}),
  \  nU_{4}(\tfrac{3}{4},  \tfrac{ \pi   n^2  }{3}),
  \ \dots,\ 
  nU_{2k}(\tfrac{3}{4},
  \tfrac{ \pi   n^2  }{3})
\end{equation}
span ${\cal{G}}_{3}(2k)$. 
Thus the equalities \eqref{key} imply the
equalities \eqref{sumGis0};
the latter are identical to \eqref{numeq}
and hence they imply 
the functional equations \eqref{thetaFE}
and \eqref{fe}.

\section{Corollaries of the functional equation
in the general case}\label{thetagen}

As was mentioned in the previous section, 
the functional identity \eqref{fe} is equivalent
to  \eqref{sumGis0}. This observation can be generalised to the case
of an arbitrary Dirichlet 
primitive character $\chi$
modulo  some~$q$, where $q>1$.

In Section \ref{thetapol} we deal with the 
character $\chi_3$, which is defined by \eqref{chi3}.
It has 
two special features:
\begin{itemize}
  \item $\chi_3$  is a real character, that is, 
  $\chi_3=\bar\chi_3$;
  \item  $\chi_3$  is an odd character, that is, 
  $\chi_3(-1)=-1$.
\end{itemize}
In the  general case 
counterparts to 
\eqref{thetaFE}--\eqref{sumGis0}
involve both $\chi$ and
$\bar\chi$, as well as $\delta$, 
the indicator of the parity of $\chi$ 
(defined by \eqref{kappa}).

More precisely,  the functional
equation takes the following form (cf.\ \eqref{thetaFE}):
\begin{equation}\label{thetaFEL}
 \theta_\chi({\tau}^{-1})=
\omega{\tau^{\delta+\frac{1}{2}}}
  \theta_{\bar{\chi}}(\tau), 
\end{equation}
where (cf.\eqref{thetaD})
\begin{equation}\label{thetaDL}
  \theta_\chi(\tau)=\sum_{n=-\infty}^\infty
  \chi(n)n^\delta\mye^{-\frac{\ \pi n  ^2}{q}\tau}
\end{equation}
and  $\omega$  is  defined by
\eqref{omega} 
(\cite[Ch.\,1, Sect. 4, Lemma\,2]{kar}).

Instead of a single series of polynomials, $G_m(n)$, 
we  
now  introduce two series:
\begin{eqnarray}\label{defE}
  E_{m}(l)&=&\mye^{{l}}\hspace{-1mm}\left.\frac{\myd^m}{\myd \tau^m} 
  \mye^{-\frac{l}{\tau}}
  \right\vert_{\tau=1}
  \\&=&\label{defE2}
  \sum_{k=0}^m (-1)^{m+k}m!k!\binom{m-1}{k-1}
  l^k,
\end{eqnarray}
\begin{eqnarray}\label{defF}
  F_{d,m}(l)&=&\mye^{l}\hspace{-1mm}\left.\frac{\myd^m}{\myd \tau^m} 
  \tau^{d+\frac{1}{2}}\mye^{-l\tau}
  \right\vert_{\tau=1}
  \\&=& \label{defF2}
    \sum_{k=0}^m(-1)^k\binom{m}{k}
 \left(d+\tfrac{1}{2}\right)^{\underline{m-k}}
 l^k,
\end{eqnarray}
where  $d=0$ or $d=1$. 
The functional equation \eqref{thetaFEL} 
is equivalent (cf.\ \eqref{numeq})
to the infinite set of 
numerical equalities
\begin{equation}\label{numeqL}
  \left.\frac{\myd^m}{\myd \tau^m} 
  \left(
 \theta_\chi({\tau}^{-1})-
\omega{\tau^{\delta+\frac{1}{2}}}
  \theta_{\bar{\chi}}(\tau)\right)
  \right\vert_{\tau=1}
  =0,\qquad m=0,1,\dots \ .
\end{equation}
In the notation \eqref{defE}--\eqref{defF},
the  equalities \eqref{numeqL} can be written as
\begin{equation}\label{sumEFL}
  \sum_{n=1}^\infty{\left(\chi(n)n ^\delta E_m\mleft(\tfrac{\pi n^2}{q}\mright)
  -\omega \bar{\chi}(n)n^\delta F_{\delta,m}\mleft(\tfrac{\pi n^2}{q}\mright)\right)\mye^{-\frac{\pi n^2}{q}}}=0.
\end{equation}
The equality \eqref{sumGis0} is a special case of 
\eqref{sumEFL}: if $\chi=\chi_3$, then $q=3$, 
$\bar\chi=\chi$, $\delta=1$, $\omega=1$, and
\begin{equation}\label{GEF}
G_m(n)=nE_m\mleft(\tfrac{\pi n^2}{3}\mright)-
nF_{1,m}\mleft(\tfrac{\pi n^2}{3}\mright).
\end{equation}

Let ${\cal{F}}_d(m) $ be the linear span of
pairs of polynomials  
\begin{equation}
  \langle E_0(l),F_{d,0}(l)\rangle, \langle E_1(l),F_{d,1}(l)\rangle,
  \dots,
  \langle E_m(l),F_{d,m}(l)\rangle.
\end{equation}

The dimension of ${\cal{F}}_d(m) $ is equal to $m+1$.

Clearly, if
$ \langle P(l),Q(l)\rangle \in {\cal{F}}_\delta(m) $
then (cf.\ \eqref{Pis0})
\begin{equation}\label{sumEF}
  \sum_{n=1}^\infty{\left(\chi(n)n^\delta  P\mleft(\tfrac{\pi n^2}{q}\mright)
  -\omega \bar{\chi}(n)n^\delta Q{}\mleft(\tfrac{\pi n^2}{q}\mright)\right)\mye^{-\frac{\pi n^2}{q}}}=0.
\end{equation}

Starting from the general representation \eqref{splitted2Low},
we can derive (by analogy with \eqref{xiU2}) that
\begin{multline}\label{xiU1gen}
  \xi_\chi(\tfrac{1}{2}+\myi t)
  =\sum_{k=1}^\infty
  \frac{2^{k}}{t^{k}}
  \sum_{n=1}^\infty \left(\chi(n)n^\delta
 {U_{k}(\tfrac{\delta}{2}+\tfrac{1}{4},
  \tfrac{ \pi   n^2  }{q})}\right.+
  \\
    \left.\ \ \ \ +(-1)^k
  \omega \bar\chi(n)n^\delta
 {U_{k}(\tfrac{\delta}{2}+\tfrac{1}{4},
  \tfrac{ \pi   n^2  }{q})}
  \right)\mye^{-\tfrac{ \pi   n^2  }{q}}.
\end{multline}

\

\textbf{Key discovery (general case).}
\emph{For $ d=0$ and for $d=1$ for every $k$}
\begin{eqnarray}\label{keygen}
\langle
 {U_{k}(\tfrac{d}{2}+\tfrac{1}{4},
  l)} ,\ -(-1)^k
 {U_{k}(\tfrac{d}{2}+\tfrac{1}{4},
 l)}
  \rangle
   \in {\cal{F}}_d(k) .
\end{eqnarray}

\ 

From \eqref{xiU1gen} and \eqref{keygen}
we obtain a generalisation of the (incorrect) equality 
\eqref{xiis0}: \emph{for every primitive Dirichlet character
$\chi$ modulo $q$ ($q>1$), for every real $t$}
\begin{equation}
\label{xiis0gen}
  \xi_\chi(\tfrac{1}{2}+\myi t)=0.
\end{equation}
Thus, the speculations from Section \ref{secplan}
can extended to an arbitrary $\chi$ as well.

\printbibliography

\begin{table}
\begin{tabular}{rl}
$U_{2}\mleft(\tfrac{3}{4},\tfrac{\pi n^2}{3}\mright)=\hspace{-3mm}$&$-\frac{\pi }{3}n^2+\frac{3}{4}$\\[2mm]
$U_{4}\mleft(\tfrac{3}{4},\tfrac{\pi n^2}{3}\mright)=\hspace{-3mm}$&$\frac{\pi ^3}{27}n^6-\frac{7 \pi ^2}{12}n^4+\frac{79 \pi }{48}n^2-\frac{27}{64}$\\[2mm]
$U_{6}\mleft(\tfrac{3}{4},\tfrac{\pi n^2}{3}\mright)=\hspace{-3mm}$&$-\frac{\pi ^5}{243}n^{10}+\frac{55 \pi ^4}{324}n^8-\frac{425 \pi ^3}{216}n^6+\frac{665 \pi ^2}{96}n^4-\frac{4141 \pi }{768}n^2+\frac{243}{1024}$\\[2mm]
$U_{8}\mleft(\tfrac{3}{4},\tfrac{\pi n^2}{3}\mright)=\hspace{-3mm}$&$\frac{\pi ^7}{2187}n^{14}-\frac{35 \pi ^6}{972}n^{12}+\frac{3689 \pi ^5}{3888}n^{10}-\frac{52745 \pi ^4}{5184}n^8+\frac{299131 \pi ^3}{6912}n^6-
$\\[1mm]&$\ \ -
\frac{185857 \pi ^2}{3072}n^4+\frac{205339 \pi }{12288}n^2-\frac{2187}{16384}$\\[2mm]
$U_{10}\mleft(\tfrac{3}{4},\tfrac{\pi n^2}{3}\mright)=\hspace{-3mm}$&$-\frac{\pi ^9}{19683}n^{18}+\frac{19 \pi ^8}{2916}n^{16}-\frac{895 \pi ^7}{2916}n^{14}+\frac{2905 \pi ^6}{432}n^{12}-\frac{744317 \pi ^5}{10368}n^{10}+
$\\[1mm]&$\ \  +\frac{4961495 \pi ^4}{13824}n^8-\frac{20501665 \pi ^3}{27648}n^6+\frac{5930365 \pi ^2}{12288}n^4-\frac{10083481 \pi }{196608}n^2+\frac{19683}{262144}$\\[2mm]
$U_{12}\mleft(\tfrac{3}{4},\tfrac{\pi n^2}{3}\mright)=\hspace{-3mm}$&$\frac{\pi ^{11}}{177147}n^{22}-\frac{253 \pi ^{10}}{236196}n^{20}+\frac{8305 \pi ^9}{104976}n^{18}-\frac{136895 \pi ^8}{46656}n^{16}+\frac{5497327 \pi ^7}{93312}n^{14}-
$\\[1mm]&$\ \ -\frac{26765585 \pi ^6}{41472}n^{12}+\frac{1847382779 \pi ^5}{497664}n^{10}-\frac{6802994495 \pi ^4}{663552}n^8+
$\\[1mm]&$\ \ +
\frac{20311855861 \pi ^3}{1769472}n^6-
\frac{2930804107 \pi ^2}{786432}n^4+\frac{494287399 \pi }{3145728}n^2-\frac{177147}{4194304}$\\[2mm]
$U_{14}\mleft(\tfrac{3}{4},\tfrac{\pi n^2}{3}\mright)=\hspace{-3mm}$&$-\frac{\pi ^{13}}{1594323}n^{26}+\frac{13 \pi ^{12}}{78732}n^{24}-\frac{24947 \pi ^{11}}{1417176}n^{22}+\frac{1871441 \pi ^{10}}{1889568}n^{20}-\frac{54017249 \pi ^9}{1679616}n^{18}+
$\\[1mm]&$\ \ +
\frac{463109933 \pi ^8}{746496}n^{16}-\frac{15833905945 \pi ^7}{2239488}n^{14}+\frac{45796856105 \pi ^6}{995328}n^{12}-
$\\[1mm]&$\ \ -
\frac{2546294373623 \pi ^5}{15925248}n^{10}+\frac{5611157536985 \pi ^4}{21233664}n^8-\frac{2399632288235 \pi ^3}{14155776}n^6+
$\\[1mm]&$\ \ +
\frac{178796540195 \pi ^2}{6291456}n^4-\frac{24221854021 \pi }{50331648}n^2+\frac{1594323}{67108864}$\\[2mm]
$U_{16}\mleft(\tfrac{3}{4},\tfrac{\pi n^2}{3}\mright)=\hspace{-3mm}$&$\frac{\pi ^{15}}{14348907}n^{30}-\frac{155 \pi ^{14}}{6377292}n^{28}+\frac{90125 \pi ^{13}}{25509168}n^{26}-\frac{117845 \pi ^{12}}{419904}n^{24}+
$\\[1mm]&$\ \ +
\frac{609126973 \pi ^{11}}{45349632}n^{22}-\frac{24376821469 \pi ^{10}}{60466176}n^{20}+\frac{618247546345 \pi ^9}{80621568}n^{18}-
$\\[1mm]&$\ \ -
\frac{364277794595 \pi ^8}{3981312}n^{16}+\frac{31762146891761 \pi ^7}{47775744}n^{14}-\frac{59284085936075 \pi ^6}{21233664}n^{12}+
$\\[1mm]&$\ \ +
\frac{529152534925469 \pi ^5}{84934656}n^{10}-\frac{729958548564245 \pi ^4}{113246208}n^8+\frac{1107751011830191 \pi ^3}{452984832}n^6-
$\\[1mm]&$\ \ -
\frac{43414094382757 \pi ^2}{201326592}n^4+\frac{1186886790259 \pi }{805306368}n^2-\frac{14348907}{1073741824}$\\[2mm]
$U_{18}\mleft(\tfrac{3}{4},\tfrac{\pi n^2}{3}\mright)=\hspace{-3mm}$&$-\frac{\pi ^{17}}{129140163}n^{34}+\frac{595 \pi ^{16}}{172186884}n^{32}-\frac{18853 \pi ^{15}}{28697814}n^{30}+\frac{297755 \pi ^{14}}{4251528}n^{28}-
$\\[1mm]&$\ \ -
\frac{157797689 \pi ^{13}}{34012224}n^{26}+\frac{1008323771 \pi ^{12}}{5038848}n^{24}-\frac{519715370239 \pi ^{11}}{90699264}n^{22}+
$\\[1mm]&$\ \ +
\frac{13200567768817 \pi ^{10}}{120932352}n^{20}-\frac{293911332044473 \pi ^9}{214990848}n^{18}+\frac{1048386106345721 \pi ^8}{95551488}n^{16}-
$\\[1mm]&$\ \ -
\frac{5171000831348105 \pi ^7}{95551488}n^{14}+\frac{6535078384876405 \pi ^6}{42467328}n^{12}-
$\\[1mm]&$\ \ -
\frac{233211900264411091 \pi ^5}{1019215872}n^{10}+
\frac{207084721195113025 \pi ^4}{1358954496}n^8-
$\\[1mm]&$\ \ -
\frac{31578932988804785 \pi ^3}{905969664}n^6+\frac{657528343229765 \pi ^2}{402653184}n^4-\frac{58157596211761 \pi }{12884901888}n^2+
$\\[1mm]&$\ \ +
\frac{129140163}{17179869184}$\\[2mm]
\end{tabular}
\caption{Polynomials $U_k(a,m)$ defined by \eqref{defU}.}
\label{U}
\end{table}

\begin{table}
\begin{tabular}{rl}
$G_{1}(n)=\hspace{-3mm}$&$\frac{2 \pi }{3}n^3-\frac{3}{2}n$\\[2mm]
$G_{2}(n)=\hspace{-3mm}$&$\frac{\pi }{3}n^3-\frac{3}{4}n$\\[2mm]
$G_{3}(n)=\hspace{-3mm}$&$\frac{2 \pi ^3}{27}n^7-\frac{7 \pi ^2}{6}n^5+\frac{11 \pi }{4}n^3+\frac{3}{8}n$\\[2mm]
$G_{4}(n)=\hspace{-3mm}$&$-\frac{2 \pi ^3}{9}n^7+\frac{7 \pi ^2}{2}n^5-\frac{17 \pi }{2}n^3-\frac{9}{16}n$\\[2mm]
$G_{5}(n)=\hspace{-3mm}$&$\frac{2 \pi ^5}{243}n^{11}-\frac{55 \pi ^4}{162}n^9+\frac{85 \pi ^3}{18}n^7-\frac{105 \pi ^2}{4}n^5+\frac{655 \pi }{16}n^3+\frac{45}{32}n$\\[2mm]
$G_{6}(n)=\hspace{-3mm}$&$-\frac{7 \pi ^5}{81}n^{11}+\frac{385 \pi ^4}{108}n^9-\frac{805 \pi ^3}{18}n^7+\frac{3185 \pi ^2}{16}n^5-\frac{3885 \pi }{16}n^3-\frac{315}{64}n$\\[2mm]
$G_{7}(n)=\hspace{-3mm}$&$\frac{2 \pi ^7}{2187}n^{15}-\frac{35 \pi ^6}{486}n^{13}+\frac{287 \pi ^5}{108}n^{11}-\frac{11165 \pi ^4}{216}n^9+\frac{22435 \pi ^3}{48}n^7-\frac{53655 \pi ^2}{32}n^5+
$\\[1mm]&$\ \ 
+\frac{108255 \pi }{64}n^3+\frac{2835}{128}n$\\[2mm]
$G_{8}(n)=\hspace{-3mm}$&$-\frac{44 \pi ^7}{2187}n^{15}+\frac{385 \pi ^6}{243}n^{13}-\frac{1309 \pi ^5}{27}n^{11}+\frac{156695 \pi ^4}{216}n^9-\frac{62755 \pi ^3}{12}n^7+\frac{250635 \pi ^2}{16}n^5
-$\\[1mm]&$\ \ 
-\frac{215985 \pi }{16}n^3-\frac{31185}{256}n$\\[2mm]
$G_{9}(n)=\hspace{-3mm}$&$\frac{2 \pi ^9}{19683}n^{19}-\frac{19 \pi ^8}{1458}n^{17}+\frac{227 \pi ^7}{243}n^{15}-\frac{6265 \pi ^6}{162}n^{13}+\frac{62741 \pi ^5}{72}n^{11}-\frac{501655 \pi ^4}{48}n^9+
$\\[1mm]&$\ \ +
\frac{1003765 \pi ^3}{16}n^7-\frac{5158125 \pi ^2}{32}n^5+\frac{31059315 \pi }{256}n^3+\frac{405405}{512}n$\\[2mm]
$G_{10}(n)=\hspace{-3mm}$&$-\frac{25 \pi ^9}{6561}n^{19}+\frac{475 \pi ^8}{972}n^{17}-\frac{6725 \pi ^7}{243}n^{15}+\frac{188125 \pi ^6}{216}n^{13}-\frac{376355 \pi ^5}{24}n^{11}+
$\\[1mm]&$\ \ +
\frac{15051575 \pi ^4}{96}n^9-\frac{12903975 \pi ^3}{16}n^7+\frac{464330475 \pi ^2}{256}n^5-\frac{310333275 \pi }{256}n^3-\frac{6081075}{1024}n$\\[2mm]
$G_{11}(n)=\hspace{-3mm}$&$\frac{2 \pi ^{11}}{177147}n^{23}-\frac{253 \pi ^{10}}{118098}n^{21}+\frac{6655 \pi ^9}{26244}n^{19}-\frac{105545 \pi ^8}{5832}n^{17}+\frac{164285 \pi ^7}{216}n^{15}-
$\\[1mm]&$\ \ -
\frac{8278655 \pi ^6}{432}n^{13}+\frac{82793095 \pi ^5}{288}n^{11}-\frac{157690225 \pi ^4}{64}n^9+\frac{2838718575 \pi ^3}{256}n^7-
$\\[1mm]&$\ \  -
\frac{11351634525 \pi ^2}{512}n^5+\frac{13647231675 \pi }{1024}n^3+\frac{103378275}{2048}n$\\[2mm]
$G_{12}(n)=\hspace{-3mm}$&$-\frac{38 \pi ^{11}}{59049}n^{23}+\frac{4807 \pi ^{10}}{39366}n^{21}-\frac{145255 \pi ^9}{13122}n^{19}+\frac{258115 \pi ^8}{432}n^{17}-\frac{2168375 \pi ^7}{108}n^{15}+
$\\[1mm]&$\ \ +
\frac{182136185 \pi ^6}{432}n^{13}-\frac{86733955 \pi ^5}{16}n^{11}+\frac{10407745425 \pi ^4}{256}n^9-\frac{20816697825 \pi ^3}{128}n^7+
$\\[1mm]&$\ \ +
\frac{149851981125 \pi ^2}{512}n^5-\frac{81852984675 \pi }{512}n^3-\frac{1964187225}{4096}n$\\[2mm]
$G_{13}(n)=\hspace{-3mm}$&$\frac{2 \pi ^{13}}{1594323}n^{27}-\frac{13 \pi ^{12}}{39366}n^{25}+\frac{767 \pi ^{11}}{13122}n^{23}-\frac{55913 \pi ^{10}}{8748}n^{21}+\frac{15101515 \pi ^9}{34992}n^{19}-
$\\[1mm]&$\ \ -
\frac{48321845 \pi ^8}{2592}n^{17}+\frac{676520845 \pi ^7}{1296}n^{15}-\frac{902016115 \pi ^6}{96}n^{13}+\frac{27060798765 \pi ^5}{256}n^{11}-
$\\[1mm]&$\ \ -
\frac{360804869425 \pi ^4}{512}n^9+\frac{1298942619975 \pi ^3}{512}n^7-\frac{4250531560275 \pi ^2}{1024}n^5+
$\\[1mm]&$\ \ +
\frac{8510470543575 \pi }{4096}n^3+\frac{41247931725}{8192}n$
\end{tabular}
\caption{Polynomials $G_m(n)$ defined by \eqref{defG}.}
\label{G}
\end{table}

\begin{table}
\begin{tabular}{rl}
$nU_{2}(\tfrac{3}{4},  \tfrac{ \pi   n^2  }{3})=\hspace{-3mm}$&$
-\frac{1}{2}G_{1}(n)$\\[4mm]
$nU_{4}(\tfrac{3}{4},  \tfrac{ \pi   n^2  }{3})=\hspace{-3mm}$&$
\frac{13}{32}G_{1}(n)+\frac{1}{2}G_{3}(n)$\\[4mm]
$nU_{6}(\tfrac{3}{4},  \tfrac{ \pi   n^2  }{3})=\hspace{-3mm}$&$
\frac{359}{512}G_{1}(n)+\frac{85}{16}G_{3}(n)-\frac{1}{2}G_{5}(n)$\\[4mm]
$nU_{8}(\tfrac{3}{4},  \tfrac{ \pi   n^2  }{3})=\hspace{-3mm}$&$
\frac{468133}{8192}G_{1}(n)+\frac{190351}{512}G_{3}(n)-\frac{1477}{32}G_{5}(n)+\frac{1}{2}G_{7}(n)$\\[4mm]
$nU_{10}(\tfrac{3}{4},  \tfrac{ \pi   n^2  }{3})=\hspace{-3mm}$&$
\frac{1881623759}{131072}G_{1}(n)+\frac{193261745}{2048}G_{3}(n)-\frac{3156027}{256}G_{5}(n)+
$\\[1mm]&$ \  \  \ +
\frac{1401}{8}G_{7}(n)-\frac{1}{2}G_{9}(n)$\\[4mm]
$nU_{12}(\tfrac{3}{4},  \tfrac{ \pi   n^2  }{3})=\hspace{-3mm}$&$
\frac{20907034805053}{2097152}G_{1}(n)+
\frac{8596877726701}{131072}G_{3}(n)-
\frac{35325623707}{4096}G_{5}(n)+
$\\[1mm]&$ \  \  \ +
\frac{32737881}{256}G_{7}(n)-
\frac{15015}{32}G_{9}(n)+
\frac{1}{2}G_{11}(n)$\\[4mm]
$nU_{14}(\tfrac{3}{4},  \tfrac{ \pi   n^2  }{3})=\hspace{-3mm}$&$
\frac{519432419559176759}{33554432}G_{1}(n)+
\frac{106803035312049895}{1048576}G_{3}(n)-
$\\[1mm]&$ \  \  \ -
\frac{1756698214958687}{131072}G_{5}(n)+
\frac{409249074377}{2048}G_{7}(n)-
\frac{390116727}{512}G_{9}(n)+
$\\[1mm]&$ \  \  \ +
\frac{16471}{16}G_{11}(n)-
\frac{1}{2}G_{13}(n)$\\[4mm]
$nU_{16}(\tfrac{3}{4},  \tfrac{ \pi   n^2  }{3})=\hspace{-3mm}$&$
\frac{25147259164376391074773}{536870912}G_{1}(n)+
\frac{10341386238841675249051}{33554432}G_{3}(n)-
$\\[1mm]&$ \  \  \ -
\frac{85053788112176057811}{2097152}G_{5}(n)+
\frac{79307042206189083}{131072}G_{7}(n)-
$\\[1mm]&$ \  \  \ -
\frac{18993638637545}{8192}G_{9}(n)+
\frac{1662115273}{512}G_{11}(n)-
\frac{63385}{32}G_{13}(n)+
\frac{1}{2}G_{15}(n)$\\[4mm]
$nU_{18}(\tfrac{3}{4},  \tfrac{ \pi   n^2  }{3})=\hspace{-3mm}$&$
\frac{2153279642725416095838675359}{8589934592}G_{1}(n)+
\frac{110687584346643190774486825}{67108864}G_{3}(n)-
$\\[1mm]&$ \  \  \ -
\frac{1820736210067848438753775}{8388608}G_{5}(n)+
\frac{848913796088552235639}{262144}G_{7}(n)-
$\\[1mm]&$ \  \  \ -
\frac{813701365343805939}{65536}G_{9}(n)+
\frac{17883551790767}{1024}G_{11}(n)-
$\\[1mm]&$ \  \  \ -
\frac{1410882207}{128}G_{13}(n)+
\frac{13889}{4}G_{15}(n)-\frac{1}{2}G_{17}(n)$\\[4mm]
$nU_{20}(\tfrac{3}{4},  \tfrac{ \pi   n^2  }{3})=\hspace{-3mm}$&$
\frac{303443163732057625427762952605293}{137438953472}G_{1}(n)+
$\\[1mm]&$ \  \  \ +
\frac{124785998078248042167134451461401}{8589934592}G_{3}(n)-
$\\[1mm]&$ \  \  \ -
\frac{256580900728526055508960344967}{134217728}G_{5}(n)+
$\\[1mm]&$ \  \  \ +
\frac{239262085449139901686517041}{8388608}G_{7}(n)-
\frac{114675835780201454785845}{1048576}G_{9}(n)+
$\\[1mm]&$ \  \  \ +
\frac{10086732814021533943}{65536}G_{11}(n)-
\frac{199809728486035}{2048}G_{13}(n)+
$\\[1mm]&$ \  \  \ +
\frac{4063696269}{128}G_{15}(n)-
\frac{181659}{32}G_{17}(n)+
\frac{1}{2}G_{19}(n)$
\end{tabular}
\caption{Polynomials $nU_{m}\mleft(\tfrac{3}{4},
  \tfrac{ \pi   n^2  }{3}\mright) $
via basis  \eqref{basisO}.}
\label{UGO}
\end{table}

\begin{table}
\begin{tabular}{rl}
$nU_{2}(\tfrac{3}{4},  \tfrac{ \pi   n^2  }{3})=\hspace{-3mm}$&$
-G_{2}(n)$\\[4mm]
$nU_{4}(\tfrac{3}{4},  \tfrac{ \pi   n^2  }{3})=\hspace{-3mm}$&$
\frac{11}{16}G_{2}(n)-\frac{1}{6}G_{4}(n)$\\[4mm]
$nU_{6}(\tfrac{3}{4},  \tfrac{ \pi   n^2  }{3})=\hspace{-3mm}$&$
-\frac{21}{256}G_{2}(n)-\frac{35}{48}G_{4}(n)+
\frac{1}{21}G_{6}(n)$\\[4mm]
$nU_{8}(\tfrac{3}{4},  \tfrac{ \pi   n^2  }{3})=\hspace{-3mm}$&$
\frac{18271}{4096}G_{2}(n)-
\frac{8337}{512}G_{4}(n)+
\frac{85}{48}G_{6}(n)-
\frac{1}{44}G_{8}(n)$\\[4mm]
$nU_{10}(\tfrac{3}{4},  \tfrac{ \pi   n^2  }{3})=\hspace{-3mm}$&$
\frac{16180951}{65536}G_{2}(n)-
\frac{6422435}{6144}G_{4}(n)+\frac{90943}{640}G_{6}(n)-
\frac{543}{176}G_{8}(n)+\frac{1}{75}G_{10}(n)$\\[4mm]
$nU_{12}(\tfrac{3}{4},  \tfrac{ \pi   n^2  }{3})=\hspace{-3mm}$&$
\frac{41997202611}{1048576}G_{2}(n)-
\frac{68523014021}{393216}G_{4}(n)+
\frac{1097088289}{43008}G_{6}(n)-
$\\[1mm]&$ \  \  \ -
\frac{342111}{512}G_{8}(n)+
\frac{1133}{240}G_{10}(n)-\frac{1}{114}G_{12}(n)$\\[4mm]
$nU_{14}(\tfrac{3}{4},  \tfrac{ \pi   n^2  }{3})=\hspace{-3mm}$&$
\frac{248669854673603}{16777216}G_{2}(n)-
\frac{68052696874435}{1048576}G_{4}(n)+
\frac{9532351392421}{983040}G_{6}(n)-
$\\[1mm]&$ \  \  \ -
\frac{1105778401}{4096}G_{8}(n)+
\frac{43310267}{19200}G_{10}(n)-
\frac{6097}{912}G_{12}(n)+\frac{1}{161}G_{14}(n)$\\[4mm]
$nU_{16}(\tfrac{3}{4},  \tfrac{ \pi   n^2  }{3})=\hspace{-3mm}$&$
\frac{2909419346984333831}{268435456}G_{2}(n)-
\frac{4783910393581167031}{100663296}G_{4}(n)+
$\\[1mm]&$ \  \  \ +
\frac{7483065768091849}{1048576}G_{6}(n)-
\frac{52982924993193}{262144}G_{8}(n)+
\frac{36494379493}{20480}G_{10}(n)-
$\\[1mm]&$ \  \  \ -
\frac{179981893}{29184}G_{12}(n)+
\frac{3305}{368}G_{14}(n)-\frac{1}{216}G_{16}(n)$\\[4mm]
$nU_{18}(\tfrac{3}{4},  \tfrac{ \pi   n^2  }{3})=\hspace{-3mm}$&$
\frac{60744738617214376444335}{4294967296}G_{2}(n)-
\frac{12489018115953506248835}{201326592}G_{4}(n)+
$\\[1mm]&$ \  \  \ +
\frac{4107272868321959300807}{440401920}G_{6}(n)-
\frac{1529918964315664377}{5767168}G_{8}(n)+
$\\[1mm]&$ \  \  \ +
\frac{5835930414241499}{2457600}G_{10}(n)-
\frac{503828721149}{58368}G_{12}(n)+
\frac{21420901}{1472}G_{14}(n)-
$\\[1mm]&$ \  \  \ -
\frac{5015}{432}G_{16}(n)+\frac{1}{279}G_{18}(n)$\\[4mm]
$nU_{20}(\tfrac{3}{4},  \tfrac{ \pi   n^2  }{3})=\hspace{-3mm}$&$
\frac{2099543857253492320406595931}{68719476736}G_{2}(n)-
\frac{1151179945136215796266665347}{8589934592}G_{4}(n)+
$\\[1mm]&$ \  \  \ +
\frac{4057413215226496747266091}{201326592}G_{6}(n)-
\frac{105904940207780519805061}{184549376}G_{8}(n)+
$\\[1mm]&$ \  \  \ +
\frac{40556259310967555519}{7864320}G_{10}(n)-
\frac{3737235414418597}{196608}G_{12}(n)+
$\\[1mm]&$ \  \  \ +
\frac{790217331055}{23552}G_{14}(n)-
\frac{709418143}{23040}G_{16}(n)+
\frac{21679}{1488}G_{18}(n)-\frac{1}{350}G_{20}(n)$\\[4mm]
\end{tabular}
\caption{Polynomials $nU_{m}\mleft(\tfrac{3}{4},
  \tfrac{ \pi   n^2  }{3}\mright) $
via basis  \eqref{basisE}.}
\label{UGE}
\end{table}

\begin{table}
\begin{tabular}{rl}
$G_{1}( n)=\hspace{-3mm}$&$-
2nU_{2}(\tfrac{3}{4},  \tfrac{ \pi   n^2  }{3})$\\[4mm]

$G_{2}( n)=\hspace{-3mm}$&$-
nU_{2}(\tfrac{3}{4},  \tfrac{ \pi   n^2  }{3})$\\[4mm]

$G_{3}( n)=\hspace{-3mm}$&$-
\frac{13}{8}nU_{2}(\tfrac{3}{4},  \tfrac{ \pi   n^2  }{3})
+2nU_{4}(\tfrac{3}{4},  \tfrac{ \pi   n^2  }{3})$\\[4mm]

$G_{4}( n)=\hspace{-3mm}$&$-
\frac{33}{8}nU_{2}(\tfrac{3}{4},  \tfrac{ \pi   n^2  }{3})-6nU_{4}(\tfrac{3}{4},  \tfrac{ \pi   n^2  }{3})$\\[4mm]

$G_{5}( n)=\hspace{-3mm}$&$-
\frac{1851}{128}nU_{2}(\tfrac{3}{4},  \tfrac{ \pi   n^2  }{3})
+\frac{85}{4}nU_{4}(\tfrac{3}{4},  \tfrac{ \pi   n^2  }{3})-2nU_{6}(\tfrac{3}{4},  \tfrac{ \pi   n^2  }{3})$\\[4mm]

$G_{6}( n)=\hspace{-3mm}$&$-
\frac{16611}{256}nU_{2}(\tfrac{3}{4},  \tfrac{ \pi   n^2  }{3})-\frac{735}{8}nU_{4}(\tfrac{3}{4},  \tfrac{ \pi   n^2  }{3})
+21nU_{6}(\tfrac{3}{4},  \tfrac{ \pi   n^2  }{3})$\\[4mm]

$G_{7}( n)=\hspace{-3mm}$&$-
\frac{727497}{2048}nU_{2}(\tfrac{3}{4},  \tfrac{ \pi   n^2  }{3})
+\frac{60739}{128}nU_{4}(\tfrac{3}{4},  \tfrac{ \pi   n^2  }{3})-\frac{1477}{8}nU_{6}(\tfrac{3}{4},  \tfrac{ \pi   n^2  }{3})
+2nU_{8}(\tfrac{3}{4},  \tfrac{ \pi   n^2  }{3})$\\[4mm]

$G_{8}( n)=\hspace{-3mm}$&$-
\frac{2351745}{1024}nU_{2}(\tfrac{3}{4},  \tfrac{ \pi   n^2  }{3})-\frac{183029}{64}nU_{4}(\tfrac{3}{4},  \tfrac{ \pi   n^2  }{3})
+\frac{6545}{4}nU_{6}(\tfrac{3}{4},  \tfrac{ \pi   n^2  }{3})-
$\\[2mm]&$\ \   \ -44nU_{8}(\tfrac{3}{4},  \tfrac{ \pi   n^2  }{3})$\\[4mm]

$G_{9}( n)=\hspace{-3mm}$&$-
\frac{561268215}{32768}nU_{2}(\tfrac{3}{4},  \tfrac{ \pi   n^2  }{3})
+\frac{10094789}{512}nU_{4}(\tfrac{3}{4},  \tfrac{ \pi   n^2  }{3})-\frac{982527}{64}nU_{6}(\tfrac{3}{4},  \tfrac{ \pi   n^2  }{3})+
$\\[2mm]&$\ \   \ 
+\frac{1401}{2}nU_{8}(\tfrac{3}{4},  \tfrac{ \pi   n^2  }{3})-2nU_{10}(\tfrac{3}{4},  \tfrac{ \pi   n^2  }{3})$\\[4mm]

$G_{10}( n)=\hspace{-3mm}$&$-
\frac{9488044485}{65536}nU_{2}(\tfrac{3}{4},  \tfrac{ \pi   n^2  }{3})-\frac{156659325}{1024}nU_{4}(\tfrac{3}{4},  \tfrac{ \pi   n^2  }{3})
+\frac{19815705}{128}nU_{6}(\tfrac{3}{4},  \tfrac{ \pi   n^2  }{3})-
$\\[2mm]&$\ \   \ -\frac{40725}{4}nU_{8}(\tfrac{3}{4},  \tfrac{ \pi   n^2  }{3})
+75nU_{10}(\tfrac{3}{4},  \tfrac{ \pi   n^2  }{3})$\\[4mm]

$G_{11}( n)=\hspace{-3mm}$&$-
\frac{717105141765}{524288}nU_{2}(\tfrac{3}{4},  \tfrac{ \pi   n^2  }{3})
+\frac{43195052901}{32768}nU_{4}(\tfrac{3}{4},  \tfrac{ \pi   n^2  }{3})-\frac{1724416375}{1024}nU_{6}(\tfrac{3}{4},  \tfrac{ \pi   n^2  }{3})+
$\\[2mm]&$\ \   \ 
+\frac{9334149}{64}nU_{8}(\tfrac{3}{4},  \tfrac{ \pi   n^2  }{3})-\frac{15015}{8}nU_{10}(\tfrac{3}{4},  \tfrac{ \pi   n^2  }{3})
+2nU_{12}(\tfrac{3}{4},  \tfrac{ \pi   n^2  }{3})$\\[4mm]

$G_{12}( n)=\hspace{-3mm}$&$-
\frac{7488983814435}{524288}nU_{2}(\tfrac{3}{4},  \tfrac{ \pi   n^2  }{3})-\frac{408626277477}{32768}nU_{4}(\tfrac{3}{4},  \tfrac{ \pi   n^2  }{3})+
$\\[2mm]&$\ \   \ 
+\frac{20218928565}{1024}nU_{6}(\tfrac{3}{4},  \tfrac{ \pi   n^2  }{3})-\frac{136171233}{64}nU_{8}(\tfrac{3}{4},  \tfrac{ \pi   n^2  }{3})+
$\\[2mm]&$\ \   \ 
+\frac{322905}{8}nU_{10}(\tfrac{3}{4},  \tfrac{ \pi   n^2  }{3})-
114nU_{12}(\tfrac{3}{4},  \tfrac{ \pi   n^2  }{3})$\\[4mm]
\end{tabular}
\caption{Polynomials $G_{m}\mleft(n\mright) $
via basis  \eqref{basisU}.}
\label{GU}
\end{table}

\end{document}